\documentclass[12pt]{article}
\usepackage{}
\usepackage{hyperref}
\usepackage{amsmath}
\usepackage{ntheorem}
\usepackage{mathtools}
\usepackage{epsf}
\usepackage{theorem}
\usepackage{indentfirst}
\usepackage{latexsym}
\usepackage{amssymb}
\usepackage[dvips]{graphicx}
\usepackage{flafter}
\usepackage{caption}
\usepackage{textcomp}
\usepackage{supertabular}
\usepackage{longtable, booktabs}
\usepackage{hhline}
\usepackage{bigstrut,bigdelim}
\usepackage{rotating}
\usepackage{multicol}
\usepackage{multirow}
\usepackage{makeidx}
\usepackage{epsfig}
\usepackage{fancyhdr}
\newif\ifger
\gerfalse
\newtheorem{theorem}{Theorem}
\newtheorem{lemma}{Lemma}[section]
\newtheorem{corollary}{Corollary}

\newtheorem{proposition}{Proposition}[section]

\begin{document}
\baselineskip=19pt

\title{Flag-transitive non-symmetric  $2$-designs with
$(r,\lambda)=1$ and exceptional groups of Lie type
}

\author{ Yongli  Zhang, Shenglin Zhou\footnote{Corresponding author. This work is
supported by the National Natural Science Foundation
of China (Grant No.11871224) and the Natural Science Foundation of Guangdong Province (Grant No.
2017A030313001). slzhou@scut.edu.cn}\\
\small \it School of Mathematics, South China University of Technology,\\
\small \it Guangzhou 510641, P.R. China}

\date{}
\maketitle

\begin{abstract}

This paper  determined all  pairs $(\mathcal{D},G)$ where $\mathcal{D}$ is a non-symmetric 2-$(v,k,\lambda)$ design   with $(r,\lambda)=1$ and  $G$ is  the  almost simple flag-transitive automorphism group of $\mathcal{D}$ with  an exceptional  socle of Lie type. We prove that if $T\trianglelefteq G\leq Aut(T)$ where $T$ is an exceptional group of Lie type, then $T$ must be the Ree group or Suzuki group, and there just five non-isomorphic designs $\mathcal{D}$.



\end{abstract}

\section{Introduction}
A 2-$(v,k,\lambda)$ design $\mathcal{D}$ is a pair $(\mathcal{P},\mathcal{B})$,
where  $\mathcal{P}$ is a set of $v$ points and  $\mathcal{B}$ is a set of $k$-subsets of $\mathcal{P}$ called blocks,
such that any 2 points are contained in exactly  $\lambda$ blocks.
 A flag is an incident  point-block pair $(\alpha, B)$.
An automorphism of   $\mathcal{D}$ is a permutation of $\mathcal{P}$  which  leaves $\mathcal{B}$ invariant.
The design is non-trivial if $2<k<v-1$ and non-symmetric if $v<b$.
All automorphisms of the  design $\mathcal{D}$ form a group called the full automorphism group of $\mathcal{D}$, denoted by $Aut(\mathcal{D})$.
Let $G\leq Aut(\mathcal{D})$, the design $\mathcal{D}$  is called  point (block, flag)-transitive
if $G$ acts transitively on the set of points (blocks, flags),  and  point-primitive if $G$ acts primitively on $\mathcal{P}$.
Note that a finite primitive group is    almost simple  if it is isomorphic to a group $G$ for which  $T\cong Inn(T)\leq G\leq Aut (T)$ for some non-abelian simple group $T$.

Let $G\leq Aut(\mathcal{D})$, and $r$ be the number of blocks incident with a given point.
In \cite{Demb1968}, P. Dembowski proved that
if $G$ is a   flag-transitive automorphism group of a  2-design $\mathcal{D}$ with $(r,\lambda)=1$, then $G$ is point-primitive.
In 1988, P. H. Zieschang \cite{Zies1988} proved that if $\cal D$ is a 2-design  with $(r,\lambda)=1$
and  $G\leq Aut(\mathcal{D})$ is flag transitive, then $G$ must be of almost simple or affine type.
Such 2-designs have been studied in \cite{MBNS, MBS, SLZZ, SLZW}, where the socle of $G$ is a sporadic, an alternating group or elementary abelian $p$-group, respectively.
In this paper, we continue to study the case that the socle of $G$ is an exceptional simple group of Lie type. We get the following:

\begin{theorem} \label{theorem}
Let $\mathcal{D}=(\mathcal{P},\mathcal{B})$ be a non-symmetric $2$-$(v,k,\lambda)$ design with $(r,\lambda)=1$
and  $G$ an  almost simple flag-transitive automorphism group of $\mathcal{D}$
with the exceptional  socle $T$ of Lie type  in characteristic $p$ and $q=p^{e}$.
Let $B$ be a block  of $\mathcal{D}$.
Then  one of the following holds:
\begin{enumerate}
\item[\rm(1)] $T={^{2}G_{2}(q)}$ with $q=3^{2n+1}\geq27$, and $\mathcal{D}$ is one of the following:
\begin{enumerate}
\item[\rm(i)] a Ree unital with $G_B=\mathbb{Z}_{2}\times L_2(q)$;
\item[\rm(ii)] a $2$-$(q^{3}+1, q, q-1)$ design with  $G_{B}=Q_{1}:K$;
\item[\rm(iii)]    a $2$-$( q^{3}+1,q, q-1)$ design with $G_{B}=Q_{2}:K$;
\item[\rm(iv)]   a $2$-$( q^{3}+1,  q^{2}, q^{2}-1)$ design with $G_{B}=Q':K$,

\end{enumerate} where  $Q\in Syl_3(T)$, and  the definitions of $Q_1, Q_2$ and $K$ refer to Section 3.
\item[\rm(2)] $T={^{2}B_{2}(q)}$ with $q=2^{2n+1}\geq8$,  and  $\mathcal{D}$ is a $2$-$(q^{2}+1, q, q-1)$ design
with  $G_{B}=Z(Q):K$, where $Q\in Syl_2(T)$ and $K=\mathbb{Z}_{q-1}\cong \mathbb{F}_{q}^*$.
\end{enumerate}
\end{theorem}

\section{Preliminary results}
We first  give some preliminary results about designs and almost simple groups.
\begin{lemma} {\rm(\cite[Lemma 2.2]{SLZZ})} \label{lemma1.2}
For a $2$-$(v,k,\lambda)$ design $\mathcal{D}$, it is well known that
\begin{enumerate}
\item[\rm(1)] \, $bk=vr$;
\item[\rm(2)] \, $\lambda(v-1)=r(k-1)$;
\item[\rm(3)] \, $v\leq \lambda v < r^{2}$;
\item[\rm(4)]\,  if $G\leq Aut(\mathcal{D})$ is flag-transitive and $(r, \lambda)=1$, then
$r\mid (|G_{\alpha}|, v-1)$ and $r\mid d $, for any non-trivial subdegree $d$ of $G$.
\end{enumerate}
\end{lemma}

\begin{lemma} \label{lemma1.1}
Assume that  $G$ and $\mathcal{D}$ satisfy the hypothesis of Theorem \ref{theorem}.
Let $\alpha\in \mathcal{P}$ and  $B\in \mathcal{B}$. Then
 \begin{enumerate}
 \item[\rm(1)]\,$G=TG_{\alpha}$ and $|G|=f|T|$ where $f$ is a divisor of $|Out(T)|$;
\item[\rm(2)] \,$|G:T|=|G_{\alpha}:T_{\alpha}|=f$;
\item[\rm(3)] \,$|G_{B}|$ divides $f|T_{B}|$, and $|G_{\alpha B}|$ divides $f|T_{\alpha B}|$ for any flag $(\alpha, B)$.
\end{enumerate}
\end{lemma}
{\bf Proof.}\,  Note that $G$ is an almost simple primitive group by $\cite{Dav1}$.  So (1) holds and (2) follows  from (1).
Since $T\unlhd G$, then $|B^{T}|$ divides $|B^{G}|$ and $|(\alpha, B)^{T}|$ divides $|(\alpha, B)^{G}|$,
hence $|G_{B}:T_{B}|$ divides $ f$,  and $|G_{\alpha B}:T_{\alpha B}|$ divides $ f$, (3) holds. $\hfill\square$

\begin{lemma} {\rm(\cite[2.2.5]{Demb1968})} \label{lemma115}
Let $\mathcal{D}$ be a $2$-$(v,k,\lambda)$ design. If $\mathcal{D}$ satisfies $r=k+\lambda$ and $\lambda\leq2$,
then $\mathcal{D}$ is embedded in a symmetric $2$-$(v+k+\lambda, k+\lambda, \lambda)$ design.
\end{lemma}
\begin{lemma} {\rm(\cite[2.3.8]{Demb1968})} \label{do}
Let $\mathcal{D}$ be a $2$-$(v,k,\lambda)$ design and $G\leq Aut(\mathcal{D})$.
If $G$ is $2$-transitive on points  and $(r,\lambda)=1$,
then $G$ is flag transitive.
\end{lemma}

\begin{lemma}  \label{lemma15}
Let $A$, $B$, $C$ be subgroups of  group $G$. If $B\leq A$, then $$|A:B|\geq|(A\cap C):(B\cap C)|.$$
\end{lemma}

\begin{lemma}{\rm(\cite{LMWS})} \label{lem:power}
 Suppose that $T$ is  a simple group of Lie type in characteristic $p$ and acts on the set of cosets of a maximal parabolic subgroup.
Then $T$ has a unique subdegree which is a power of $p$ except $T$ is $ L_d(q)$, $\Omega^{+}_{2m}(q)$ ($m$ is odd) or $E_6(q)$.
\end{lemma}

\begin{lemma}{\rm{\cite[1.6]{GMS}}}$\mathrm{(Tits~Lemma)}$\label{Tits}
If $T$ is  a simple group of Lie type in characteristic $p$, then any proper subgroup of index prime to $p$ is contained in a  parabolic subgroup of $T$.
\end{lemma}

In the following, for a positive integer $n$, $n_p$ denotes the $p$-part of $n$ and $n_{p'}$ denotes the $p'$-part of $n$, i.e., $n_p=p^t$ where $p^t\mid n$ but $p^{t+1}\nmid n$, and $n_{p'}=n/n_{p}$.
\begin{lemma}  \label{lifanggen}
Assume that  $G$ and $\mathcal{D}$ satisfy the hypothesis of Theorem \ref{theorem}. Then $|G|<|G_\alpha|^3$ and if $G_{\alpha}$ is   a non-parabolic maximal subgroup of $G$, then
$|G|<|G_\alpha|{|G_\alpha|_{p'}^{2}}$ and $|T|<|Out(T)|^2|T_\alpha|{|T_\alpha|_{p'}^{2}}$.
\end{lemma}
{\bf Proof.}\, From Lemma  \ref{lemma1.2}, since $r$ divides every non-trivial subdegree of $G$,
then $r $ divides $|G_{\alpha}|$, and so   $|G|<|G_\alpha|^3$.
If $G_{\alpha}$ is  not  parabolic, then   $p$ divides $v=|G:G_\alpha|$ by Lemma \ref{Tits}.
Since $r$ divides $v-1$, $(r,p)=1$ and so $r$ divides $|G_\alpha|_{p'}$.
It follows that $r<|G_\alpha|_{p'}$, and hence  $|G|<|G_\alpha|{|G_\alpha|_{p'}^{2}}$ by Lemma \ref{lemma1.2}.
Now by  Lemma \ref{lemma1.1}(2), we have that  $|T|<|Out(T)|^2|T_\alpha|{|T_\alpha|_{p'}^{2}}$.
$\hfill\square$

\begin{lemma} {\rm(\cite[Theorem 2, Table III]{Lie90})} \label{theorem 10}
If $T$ is a finite simple exceptional group of Lie type such that $T\leq G \leq Aut(T)$, and $G_\alpha$ is a maximal subgroup of $G$ such that $T_0=Soc(G_\alpha)$ is not simple, then one of the following holds:
\begin{enumerate}
\item[\rm(1)] \, $G_\alpha$ is parabolic;
\item[\rm(2)] \,$G_\alpha$ is of maximal rank;
\item[\rm(3)] \, $G_\alpha= N_G(E)$, where $E$ is an elementary abelian group given in \cite[Theorem 1 (II)]{CLSS92};
\item[\rm(4)] \, $T=E_8(q)$ with $p>5$, and $T_0$ is either $A_5\times A_6$ or $A_5\times L_2(q)$;
\item[\rm(5)] \, $T_0$ is as in Table \ref{Table 1}.
\end{enumerate}
\begin{table}[ht]\caption{ }
\label{Table 1}
\noindent\[
\begin{array}{ll}
\hline
\mbox{$T$}&\mbox{\hspace{4cm}$T_0$}\\
\hline
F_4(q)&L_2(q)\times G_2(q) (p>2, q>3)\\
E_6^{\epsilon}(q)&L_3(q)\times G_2(q), U_3(q)\times G_2(q) (q>2)\\
E_7(q)&L_2(q)\times L_2(q) (p>3), L_2(q)\times G_2(q) (p>2,q>3),\\
     &L_2(q)\times F_4(q) (q>3), G_2(q)\times Sp_6(q)\\
E_8(q)& L_2(q)\times L_3^\epsilon(q) (p>3), L_2(q)\times
G_2(q)\times G_2(q) (p>2, q>3),\\
      &G_2(q)\times F_4(q), L_2(q)\times
      G_2(q^2) (p>2, q>3)\\
\hline
\end{array}
\]
\end{table}
\end{lemma}

\begin{lemma} {\rm(\cite[Theorem 3]{LST96})} \label{theorem 11}
Let $T$ be a finite simple exceptional group of Lie type, with $T\leq G \leq Aut(T)$. Assume  $G_\alpha$ is a maximal subgroup of $G$ and $ Soc(G_\alpha)=T_0(q)$ is a simple group of Lie type over $\mathbb{F}_q (q>2)$ such that $\frac{1}{2}rank(T)<rank(T_0)$; assume also that $(T, T_0)$ is not $(E_8,{^2A_5(5)})$ or $(E_8,{^2D_5(3)})$.
Then one of  the following holds:
\begin{enumerate}
\item[\rm(1)] \, $G_\alpha$ is a subgroup of maximal rank;
\item[\rm(2)] \,$T_0$ is  a subfield or twisted subgroup;
\item[\rm(3)] \, $T=E_6(q)$ and $T_0=C_4(q) (q~odd)$ or $F_4(q)$.
\end{enumerate}
\end{lemma}

\begin{lemma} {\rm(\cite[Theorem 1.2]{LS95})} \label{theorem 12}
Let $T$ be a finite simple exceptional group of Lie type such that $T\leq G \leq Aut(T)$, and $G_\alpha$  a maximal subgroup of $G$ with socle  $T_0=T_0(q)$ a simple group of Lie type in characteristic $p$.
Then  if $rank(T_0)\leq \frac{1}{2}rank(T)$, we have the following bounds:
\begin{enumerate}
\item[\rm(1)] \,if $T=F_4(q)$, then $|G_\alpha|<4 q^{20}\log_pq ;$
\item[\rm(2)] \,if $T=E_{6}^{\epsilon}(q)$, then $|G_\alpha|<4 q^{28}\log_pq$;
\item[\rm(3)] \,if $T=E_7(q)$, then $|G_\alpha|<4 q^{30}\log_pq$;
\item[\rm(4)] \,if $T=E_8(q)$, then $|G_\alpha|<12 q^{56}\log_pq$.
\end{enumerate}
In all cases, $|G_\alpha|<12 |G|^{\frac{5}{13}}\log_pq.$
\end{lemma}

The following lemma  gives a method  to check the existence of the  design with  the  possible parameters.
\begin{lemma} \label{lemma2.1}
For the given parameters $(v,b,r,k,\lambda)$ and the  group $G$, the conditions that there exists a design $\mathcal{D}$ with such parameters
satisfying $G$  which is flag-transitive and point primitive is equivalent to the following four steps holding
for some subgroup $H$ of  $G$ with  index $b$ and its orbit of size $k$:
\begin{enumerate}
\item[\rm(1)]  $G$ has at least one subgroup $H$ of order $|G|/b$;
\item[\rm(2)]  $H$ has at least one orbit $O$ of length $k$;
\item[\rm(3)]  the size of $O^{G}$ is $b$;
\item[\rm(4)]  the number of blocks incident with any two points  is a constant.
\end{enumerate}

\end{lemma}

When we run through all possibilities of $H$ and its  orbits  with size $k$, then we found all designs with such parameters
and admitting $G\leq Aut (\mathcal{D})$  is flag-transitive and point primitive.
This is the essentially strategy adopted in \cite{SLZZ}.

We now give some information about the  Ree group $^{2}G_{2}(q)$ with $q=3^{2n+1}$ and  its subgroups,
which from \cite{XGFLR,BH1982,SREE} and would be used later.

Set $m=3^{n+1}$, and so $m^{2}=3q$. The  Ree group $^{2}G_{2}(q)$  is   generated by $Q, K$ and $\tau$,
where $Q$ is Sylow 3-subgroup of $^{2}G_{2}(q)$,
$K=\{diag(t^{m},t^{1-m},t^{2m-1},1,t^{1-2m},t^{m-1}, t^{-m})\,|\, t\in\mathbb{F}_{q}^{*}\}\cong \mathbb{Z}_{q-1}$
and $\tau^{2}=1$ such that $\tau$ inverts $K$, and $|^{2}G_{2}(q)|=(q^{3}+1)q^{3}(q-1)$.
\begin{lemma} \label{rpro}
\begin{enumerate}
\item[\rm(1)] \,{\rm(\cite{SREE})}   $^{2}G_{2}(q)$  is $2$-transitive of degree $q^{3}+1$.
\item[\rm(2)] \,{\rm(\cite[p.252]{pg})} The stabilizer of one point  is $Q:K$, and $N_{^{2}G_{2}(q)}(Q)=Q:K$.
\item[\rm(3)] \,{\rm(\cite[p.292]{BH1982})} The stabilizer $K$  of two points is cyclic of order $q-1$ and the stabilizer of three points is of order$2$.
\item[\rm(4)] \,{\rm(\cite[p.292]{BH1982})} The Sylow $2$-subgroup of $^{2}G_{2}(q)$ is   elementary abelian with order $8$.

\end{enumerate}
\end{lemma}

\begin{lemma}{\rm(\cite[Lemma 3.3]{XGFLR} )} \label{REW}
Let $M\leq {^{2}G_{2}(q)}$ and $M$ be maximal in $^{2}G_{2}(q)$.
Then either $M$ is conjugate to $M_{6}:={^{2}G_{2}(3^{\ell})}$ for some divisor $\ell$ of $2n+1$, or
$M$ is conjugate to one of the subgroups $M_{i}$ in the following table:
  \begin{table}[htbp]
 \centering
\caption{ The maximal subgroups of $^{2}G_{2}(q)$}\label{rank}

\begin{tabular}{cccc}
\toprule
Group&~~Structure & ~~Remarks\\
\midrule
$M_{1}$&~~$Q:K$&~~the normalizer of $Q$ in $^{2}G_{2}(q)$\\
$M_{2}$&$\mathbb{Z}_{2}\times L_{2}(q)$&~~the centralizer of an involution in $^{2}G_{2}(q)$\\
$M_{3}$&$(\mathbb{Z}_{2}^{2}\times D_{(q+1)/2}):\mathbb{Z}_{3}$&~~the normalizer of a four-subgroup\\
$M_{4}$&$\mathbb{Z}_{q+m+1}:\mathbb{Z}_{6}$&~~the normalizer of $\mathbb{Z}_{q+m+1}$ \\
$M_{5}$&$\mathbb{Z}_{q-m+1}:\mathbb{Z}_{6}$&~~the normalizer of $\mathbb{Z}_{q-m+1}$ \\
\bottomrule
\end{tabular}
\end{table}

\end{lemma}

Moreover, we see that from  \cite{XGFLR}, the Sylow 3-subgroup $Q$ can
be identified with the group consisting of all triples $(\alpha,\beta,\gamma)$
from  $\mathbb{F}_{q}$ with   multiplication:
$$(\alpha_{1},\beta_{1}, \gamma_{1})(\alpha_{2}, \beta_{2}, \gamma_{2})
= (\alpha_{1}+\alpha_{2}, \beta_{1}+\beta_{2}-\alpha_{1}\alpha^{m}_{2}, \gamma_{1}+\gamma_{2}- \alpha^{m}_{1}\alpha^{m}_{2}-\alpha_{2}\beta_{1}+ \alpha_{1}\alpha^{m+1}_{2}).$$ It is easy to check that $(0,0,\gamma)(0,\beta,0)=(0,\beta,\gamma)$.
Set $Q_{1}=\{(0,0,\gamma)|\gamma\in \mathbb{F}_{q}\}$ and
$Q_{2}=\{(0,\beta,0)|\beta\in \mathbb{F}_{q}\}$,
then $Q_{1}\cong Q_{2}\cong \mathbb{Z}_{3}^{2n+1}$.

For a group $Q$, $Z(Q)$, $\Phi(Q)$, $Q'$ denote the center, Frattini subgroup, and the derived subgroup  of $Q$, respectively.
Then  $Q'=\Phi(Q)=Q_{1}\times Q_{2}$, $Z(Q)=Q_{1}$, and $Q'$ is an elementary abelian $3$-group.
For any $(\alpha,\beta,\gamma)\in Q$ and $k\in K$,
$$(\alpha,\beta,\gamma)^{k}=(k\alpha,~ k^{1+m}\beta,~ k^{2+m}\gamma).$$

\begin{lemma}{\rm(\cite{XGFLR,SREE})} \label{REQ}
Let $Q$, $M$, $Q_2$, $M_2$ and $K$ as above, then
\begin{enumerate}
\item[\rm(1)] the normalizer of  any subgroup of $Q$ is contained in $M_{1}$;
\item[\rm(2)] for any $g\in {^{2}G_{2}(q)}$, either $Q^g=Q$ or $Q^g\cap Q=1$;
\item[\rm(3)] $Q_{2}$ is a Sylow $3$-subgroup of $M_{2}$  and $N_{M_{2} }(Q_{2})=2\times (Q_{2}:\langle k^{2}\rangle)$ with $\langle k\rangle=K$.
\end{enumerate}
\end{lemma}

\begin{lemma}{\rm(\cite[Lemma 3.2]{XGFLR})} \label{REK}
The following  hold for  the cyclic subgroup $K$:
\begin{enumerate}
\item[\rm(1)]  $K$ is transitive  on $Q_{1}\setminus \{1\}$ acting by conjugation;
\item[\rm(2)]  $K$ has two orbits $(0,1,0)^{K}$, $(0,-1,0)^{K}$ on $Q_{2}\setminus \{1\}$ acting by conjugation.
\end{enumerate}
\end{lemma}

From above lemmas, we have the following properties of the subgroups of $^{2}G_{2}(q)$.
\begin{lemma} \label{AK5}
If $H\leq M_{1}$ and  $(q-1)\mid |H|$, then   $K\leq H$.
\end{lemma}
{\bf Proof.}\, Let $p$ be a prime divisor of $ q-1$.  If $P\in Syl_{p}(M_{1})$,
then since $(p,3)=1$ and $Q\cap K=1$, we have $P\in Syl_{p}(K)$.
Note that $K$ is  cyclic,  the Sylow $p$-subgroup of $K$ is unique,
and so the Sylow $p$-subgroup of $M_{1}$ is unique.
On the other hand, if $P_0\in Syl_{p}(H)$, since $H\leq M_{1}$, then $P_{0}=P\cap H$.
Moreover, $|P_{0}|=|P|$ implies that $P=P_{0}\leq H$.
Since  $p$ is arbitrary, all Sylow subgroups of $K$  are contained in $H$, and so   $K\leq H$.
$\hfill\square$

\begin{corollary} \label{AK1}
Let $H\leq M_{1}$ and  $|H|=q(q-1)$. Then $H=A:K$ where $A$ is the Sylow $3$-subgroup of $H$.
\end{corollary}
{\bf Proof.} Since  $Q\unlhd M_{1}$, we have $A= H\cap Q$ and $A\unlhd H$.
By Lemma \ref{AK5}, $K\leq H$. Now $A \cap K=1$, and so $H=A:K$.
$\hfill\square$

\begin{lemma} \label{AK14} Let $Q_{2}$ be a Sylow $3$-subgroup of $M_{2}$ and $H_{2}:=N_{M_{2}}(Q_{2}) $.
If $Q_{2}\leq Q$  and $M_{1}=Q:K$, then the following hold:
\begin{enumerate}
\item[\rm(1)]  $H_{2}=Q_{2}:K$ and  $H_{2}\leq M_{1}$;
\item[\rm(2)] for any $H\leq M_{2}$ satisfying   $|H|=q(q-1)$, there exists $c\in M_{2}$ such that $H=H_{2}^{c}$ and $H\leq M_{1}^{c}$ .
\end{enumerate}

\end{lemma}
{\bf Proof.}\, Clearly, (1) holds by Lemma \ref{REQ}(1) and Corollary \ref{AK1}.
Let $H\leq M_{2}$ and $|H|=q(q-1)$. Note that  $M_{2}\cong\mathbb{Z}_{2}\times  L_{2}(q)$.
Since $H\lesssim \mathbb{Z}_{2}\times  L_{2}(q)$ and  $H_{2}\lesssim\mathbb{Z}_{2}\times  L_{2}(q)$,
then by the list of  maximal subgroups of $ L_{2}(q)$, we know that  $H\cong H_{2}\cong \mathbb{Z}_{2}\times ([q]:Z_{\frac{q-1}{2}})$.
Let $\sigma$ be an automorphism from $H_2$ to $H$. Then $Q_{2}^\sigma \unlhd H$ since $ Q_{2}\unlhd H_{2}$.
Moreover, since  $q\mid |H|$, the Sylow 3-subgroup of $H$ is conjugate to $Q_{2}$ in $M_{2}$
and so $Q_{2}^{\sigma}=Q_{2}^{c}\unlhd H$ for $c\in M_{2}$.
It follows that $$H\leq N_{M_{2}}(Q_{2}^{c})=N_{M_{2}}(Q_{2})^{c}=H_{2}^{c}.$$
Therefore  $H=H_{2}^{c}$.

Note that if $Q^{c}\neq Q$, then from $Q_{2}^{c}\leq Q^{c}$ and Lemma \ref{REQ}(1),
we get $H= N_{M_{2}}(Q_{2}^{c})\leq M_{1}^{c}$, and so (2) holds.

Now, we prove that $Q^{c}\neq Q$. If $Q^{c}= Q$, then $Q_{2}^{c}\leq Q$, and so $H\leq M_{1}$.
By Corollary \ref{AK1}, we have $H=Q_{2}^{c}:K$ and $H_{2}=Q_{2}:K$.
Since $Q_{2}\unlhd Q'$, $Q_{2}^{c}\unlhd Q'$.  Recall that  $Q'=Q_1\times Q_2$ is an elementary abelian 3-group,
so $ Q_{2}^{c}\cap Q_{1}\neq1$  or $ Q_{2}^{c}\cap Q_{2}\neq1$.
Now suppose that $(0,\beta,0)\in Q_{2}^{c}\cap Q_{2}$, since $Q_{2}^{c}\cap Q_{2}\leq Q_{2}$,
we have $(0,\beta,0)^{-1}=(0,-\beta,0)\in  Q_{2}^{c}\cap Q_{2}$.
This, together with $K\leq H $ and $K\leq H_{2} $,
implies $(0,\beta,0)^{K}\cup (0,-\beta,0)^{K}=Q_{2}\setminus \{1\}= Q_{2}^{c}\setminus \{1\}$.
Hence  $Q_{2}^{c}=Q_{2}$,  a contradiction.
Similarly, if $ Q_{2}^{c}\cap Q_{1}\neq1$, we have $ Q_{2}^{c}=Q_{1}$,  a contradiction. $\hfill\square$

\begin{lemma} \label{AK2}
Suppose that $H\leq {^{2}G_{2}(q)}$ and  $|H|=q(q-1)$. Then $H$ is conjugate to $H_{1}=Q_{1}:K$ or $H_{2}=Q_{2}:K$,
and there are only two  conjugacy classes of subgroups  of order $q(q-1)$ in $^{2}G_{2}(q)$.
\end{lemma}
{\bf Proof.}\,  Let $H\leq {^{2}G_{2}(q)}$ and $|H|=q(q-1)$.
By Lemma \ref{REW}, $H$ must be contained in a conjugacy of $M_1$ or $M_2$.
Firstly, if  $H^{g^{-1}} \leq M_1$, then by Corollary \ref{AK1},  $H^{g^{-1}}=A:K$ where $A$ is a Sylow 3-subgroup of $H^{g^{-1}}$.
We now show that  $A\leq Q' $. Assume that $F$ is a   maximal subgroup of $Q$ such that $A\leq F$.
If  $A\cap Q'=1$, then by Lemma \ref{lemma15} and the fact $Q'\leq F$, we have $|F:A|\geq |F\cap Q':A\cap Q'|= q^{2},$
and so $|F|\geq q^{3}$, a contradiction.
Therefore, there exists an element  $(0,\beta, \gamma)\in A\cap Q'$,
which implies that $A\setminus \{1\}=(0,\beta, \gamma)^{K}\subseteq Q'\setminus \{1\}$ and hence $A\leq Q'$.
It follows that  $A\cap Q_{1}\neq1$ or $A\cap Q_{2}\neq1$.
Similar to the proof of Lemma \ref{AK14}, if $A\cap Q_{1}\neq1$, then $A= Q_{1}$ and so $H^{g^{-1}}=H_1$,
and if  $A\cap Q_{2}\neq1$, then $A= Q_{2}$ and so $H^{g^{-1}}=H_2$.
Secondly, if  $H$ contained in a conjugacy of $M_2$, then  $H$ is conjugate to $H_2$ by  Lemma \ref{AK14}(2).$\hfill\square$

\begin{lemma} \label{AK3}
Let $H\leq{^{2}G_{2}(q)}$ and  $|H|=q^{2}(q-1)$. Then $H$ is conjugate to $Q':K$,  and there are only one conjugacy class of subgroups of order $q^{2}(q-1)$ in $^{2}G_{2}(q)$.
\end{lemma}
{\bf Proof.}\, Since $Q'~char~Q\unlhd M_{1}$, so $Q':K$ is a subgroup of $ M_{1}$  with order $q^{2}(q-1)$.
Suppose that $H\leq {^{2}G_{2}(q)}$ and $|H|=q^{2}(q-1)$. By Lemma \ref{REW},  we have $H^{g^{-1}} \leq M_1$.
Similarly as the proof of Corollary \ref{AK1}, we get that $H^{g^{-1}}$ has the structure $A:K$ where $A$ is the Sylow 3-subgroup of $H^{g^{-1}}$.
Let  $F$ be a maximal subgroup of $Q$ satisfying  $A\leq F$.
Since $|F:A|\geq |F\cap Q_{i}:A\cap Q_{i}|$, we have $|A\cap Q_{i}|>1$, which implies $Q_{i}=Q_{i}^{K}\leq A^{K}=A$ for $i=1,2$.
So $Q'\leq A$, and it follows that $Q'=A$ and $H^{g^{-1}} = Q':K$ in $M_{1}$.
$\hfill\square$

Similarly, we have the following result on the Suzuki group $^{2}B_{2}(q)$ by \cite{XGFS} and \cite[p.250]{pg}.

\begin{lemma} \label{AK4}
Suppose that  $Q$ is the Sylow $2$-subgroup of  $^{2}B_{2}(q)$ and  $M_{1}=Q:K$ is the normalizer of $Q$.
Let $H\leq {^{2}B_{2}(q)}$ and $|H|=q(q-1)$.  Then $H$ is conjugate to $Z(Q):K$.
There exists a unique conjugacy class of subgroups of order $q(q-1)$ in $^{2}B_{2}(q)$.
\end{lemma}

\section{Proof of Theorem 1}

\subsection{T is  the Ree group}

\begin{proposition} \label{pro1}
Suppose that $G$ and $\mathcal{D}$ satisfy hypothesis of Theorem \ref{theorem}. Let $B$ be a block.
 If $T={}^{2}G_{2}(q)$  with $q=3^{2n+1}$,
then  $\cal D$ is the  Ree unital  or one of the following:
\begin{enumerate}
\item[\rm(1)]  $\mathcal{D}$ is a $2$-$(q^{3}+1, q, q-1)$ design with  $G_B=Q_{1}:K$ or $Q_{2}:K$;
\item[\rm(2)]  $\mathcal{D}$ is a $2$-$( q^{3}+1,  q^{2}, q^{2}-1)$ with $G_B=Q':K$.
\end{enumerate}
\end{proposition}

This proposition will be proved  into two steps. We first assume that there exists a design satisfying the assumptions
and   obtain the possible  parameters  $(v, b, r, k, \lambda)$  in Lemma \ref{lmmm1},
then  prove the existence of the  designs using Lemma \ref{lemma2.1}.

\begin{lemma} \label{lmmm1}
Suppose that $G$ and $\mathcal{D}$ satisfy the hypothesis of  Theorem \ref{theorem}. If $T={}^{2}G_{2}(q)$ with $q=3^{2n+1}$,
then $(v, b, r,k,\lambda)=( q^{3}+1,q^{2}(q^{3}+1), q^{3}, q, q-1)$ or $( q^{3}+1, q(q^{3}+1),q^{3}, q^{2}, q^{2}-1)$ or $\cal D$ is the  Ree unital.
\end{lemma}
{\bf Proof.}\, Let $T_{\alpha}:=G_\alpha \cap T$. Since $G$  is primitive on $\mathcal{P}$, then $T_\alpha$ is one of the cases in Lemma \ref{REW} by \cite{PBK1}.
First, the cases that   $T_\alpha=\mathbb{Z}_{2}^{2}\times D_{(q+1)/2}$ and  $\mathbb{Z}_{q\pm m+1}:\mathbb{Z}_{6}$
are impossible by Lemma \ref{lifanggen}.
If  $T_\alpha =\mathbb{Z}_{2}\times L_{2}(q)$, then $v=q^2(q^2-q+1)$ and $(|G_\alpha \cap T|, v-1)=(q(q^2-1),q^4-q^3+q^2-1)=q-1$.
But  since  $r$ divides $f(|G_\alpha \cap T|, v-1)$,  which is too small to satisfy $v<r^2$.
Similarly,   $T_\alpha $ cannot be ${^{2}G_{2}(3^{\ell})}$.

We next assume that $T_\alpha =Q:K$, and so $v=q^{3}+1$.
Moreover, from \cite[p.252]{pg}, $T$ is 2-transitive on $\mathcal{P}$, so $T$ is flag-transitive by Lemma \ref{do}.
Hence we may assume that $G=T={}^{2}G_{2}(q)$.
The equations in Lemma \ref{lemma1.2}  show $$b=\frac{\lambda v(v-1)}{k (k -1)}=\frac{\lambda q ^3(q^3+1)}{k(k -1)},$$
then by the flag-transitivity of $T$, we have $$|T_{B}|=\frac{|T|}{b}=\frac{(q-1)k(k-1)}{\lambda}.$$
Let $M$ be a maximal subgroup of $T$ such that $T_{B}\leq M$. Then since  $|T_{B}|  \mid   |M|$ and $q\geq27$,
$M$ must be  $M_{1}$ or  $M_{2}$ shown in Lemma \ref{REW}.

If   $T_{B}\leq M_{1} $, then   $k(k-1)\mid\lambda q^3$.
Furthermore, since $(r,\lambda)=1$ and so $\lambda\mid (k-1)$ by  Lemma \ref{lemma1.2}(2).
Therefore  $\lambda=k-1$, and it follows that $r=v-1=q^{3}$ and  $k\mid q^{3}$.
Note that $M_{1}$ is   point stabilizer of $T$ in this  action.  So there exists $\alpha$
such that $M_{1}=T_{\alpha}$ and $T_{B}\leq T_{\alpha}$.
However, the  flag-transitivity of $T$ implies   $\alpha\notin B$.
For any point $\gamma \in B$, $T_{\gamma B}\leq T_{\alpha\gamma}$.
By Lemma \ref{rpro}, $|T_{\alpha\gamma}|=q-1$,
 and so  $|T_{\gamma B}|\mid (q-1)$.
On the other hand, from $$|B^{T_{\gamma}}|=|T_{\gamma}:T_{\gamma B}|\leq |B^{G_{\gamma}}|=|G_{\gamma}:G_{\gamma B}|=r=q^{3},$$
we have  $T_{\gamma B}=T_{\alpha\gamma}$ and so $B^{T_{\alpha\gamma}}=B$.
Since the stabilizer of three points is of order 2 by Lemma \ref{rpro},
so the size of $T_{\alpha\gamma}$-orbits  acting on $\mathcal{P} \setminus \{\alpha,\gamma\}$ is $q-1$ or $\frac{1}{2}(q-1)$. This , together with  $B^{T_{\alpha\gamma}}=B$ and $\alpha\notin B$,
implies that  $k-1=a\frac{(q-1)}{2}$ for an integer $a$.
Recall that $k\mid q^{3}$ and $k< r$, we get $k=q$ or $k=q^{2}$.
If $k=q$,  then  $$b=q^{2}(q^{3}+1), r=q^{3}, \lambda=q-1.$$
If $k=q^{2}$, we have $$b=q(q^{3}+1), r=q^{3}, \lambda=q^{2}-1.$$

Now we deal with the case   that $T_{B}\leq M_{2}$ by the similar method in \cite[Theorem 3.2]{PBKF}.

If $T_{B}$ is a solvable subgroup of $M_{2} \cong \mathbb{Z}_{2}\times L_{2}(q) $,
then  $T _{B}$ must map into either  $\mathbb{Z}_{2}\times A_{4}$, $\mathbb{Z}_{2}\times D_{q\pm 1}$ or $\mathbb{Z}_{2}\times ([q]:\mathbb{Z}_{\frac{q-1}{2}})$.
Obviously, the former two cases are impossible.
For the last case, $T _{B}\lesssim \mathbb{Z}_{2}\times([q]:\mathbb{Z}_{\frac{q-1}{2}})$.
Since  $T _{B}\leq M_{2}$, by Lemma \ref{AK14}, this   can be reduced to the case  $ T _{B}\leq M_{1}$.

If $T_{B}$ is non-solvable, then it embeds in $\mathbb{Z}_{2}\times  L_{2}(q_{0})$ with $q_{0}^{\ell}=q=3^{2n+1}$.
The condition  that $|T_{B}|$ divides $ |\mathbb{Z}_{2}\times  L_{2}(q_{0})|$ forces $q_0=q$ and  so $T_{B}$ is isomorphic to $\mathbb{Z}_{2}\times  L_{2}(q)$ or  $ L_{2}(q)$.

If $T_{B} \cong\mathbb{Z}_{2}\times  L_{2}(q)$, then $T_{B}=M_{2}$ and so $b=q^{2}(q^{2}-q+1)$. Hence, from Lemma \ref{lemma1.2},
we have $k\mid q(q+1),~ q^{2}\mid r$ and  $r\mid q^{3}.$
Since  $k\geq3$, then the fact that the stabilizer of three points is of order 2 implies that  $T_{B}$ cannot acting trivially on the block $B$.
Moreover, since  $q+1$ is the smallest degree of any non-trivial action of $ L_{2}(q)$, we have  $k=\frac{\lambda(v-1)}{r}+1\geq q+1$.

If the design $\mathcal{D}$ is a linear space, then $\mathcal{D}$ is the Ree unital (see \cite{PBKF})
with parameters  $$(v, b, r,k,\lambda)=( q^{3}+1,q^{2}(q^{2}-q+1), q^{2}, q+1, 1)$$
and $T$ is flag-transitive  with the block stabilizer $M_{2}$.

If $\lambda>1$, we claim that $\lambda=k-1$.
Clearly, $\lambda\mid (k-1)$ as $(r,\lambda)=1$ by Lemma \ref{lemma1.2}(2).
If $3\mid (k-1)$ and $(k,3)=1$, then since $k\mid q(q+1)$ and  $k\geq q+1$, we have $k=q+1$ and so $\lambda \mid q$,
which  contradicts   $(r,\lambda)=1$ as $q^{2}\mid r$. Hence $(k-1,3)=1$.
Moreover, $(k-1)\mid \lambda q^{3}$ implies that $(k-1)\mid \lambda$. So we have $\lambda=k-1$.

Let $\Delta_{1}$, $\Delta_{2}$,$\ldots$, $\Delta_{t}$ be the orbits of $M_{2}$.
Since $M_{2}$ is the block stabilizer of the Ree unital, it has an orbit  of size $q+1$.
Without loss of generality, suppose that  $|\Delta_{1}|=q+1$.
On the one hand, recall that $k\mid q(q+1)$ and  $T$ is flag transitive, $T_{B}=M_2$ has at least one orbit with   size    less than $q(q+1)$.
On the other hand, we show that $|\Delta_{i}|>q(q+1)$ for $i\neq1$ in the following and we obtain the desired contradiction.
Assume that $\delta\in \mathcal{P}\setminus\Delta_{1}$, we claim that $({M_{2}})_{\delta}$ is a 2-group.
Let $p $ be a prime divisor of $ |({M_{2}})_{\delta}|$  and $P$ be a Sylow $p$-subgroup of $({M_{2}})_{\delta}$.
If $p\neq2$ and $p\neq3$, then since $({M_{2}})_{\delta}\leq T_{\delta}$,  we have $p\mid (q-1)$.
Obviously, since $\Delta_{1}$ is an orbit of $M_{2}$ and $P\leq({M_{2}})_{\delta}$, and so $P$  acts invariantly on $\Delta_{1}$ and $\mathcal{P}\setminus\Delta_{1}$.
Note that  the length of a $P$-orbit is either 1 or divided by $p$, so $P$ fixes at least two points in $\Delta_{1}$.
Moreover, $P$ also fixes $\delta$. Therefore  $P$ fixes at least three points of $\mathcal{P}$,
which is impossible as the order of the stabilizer of three points  is 2 by Lemma \ref{rpro}(3).
If $p=3$, since  $P$  fixes the point $\delta\in \mathcal{P}\setminus\Delta_{1}$ and $|\mathcal{P}\setminus\Delta_{1}|=q^{3}-q$,
then $P$ fixes at least three points in $\mathcal{P}\setminus\Delta_{1}$, which is also impossible.
As a result, $({M_{2}})_{\delta}$ is a 2-group. The fact that  the Sylow 2-subgroup of $T$ is of  order 8 implies that
 the sizes of the $M_2$-orbits $\Delta_{i}$ ($i\neq 1$) are at least $\frac{q(q^{2}-1)}{8}$ and hence larger than $q(q+1)$,
which contradicts  the fact $k\mid q(q+1)$. Therefore, $T_{B} \not\cong\mathbb{Z}_{2}\times  L_{2}(q)$. Similarly, $T_{B}\not \cong L_{2}(q)$.
Thus  $T_{B}$ is not a  non-solvable subgroup in $M_2$. $\hfill\square$

{\bf Proof of Proposition 3.1.}\ We   use    Lemma \ref{lemma2.1} to prove  the existence of the design  with parameters listed in Lemma \ref{lmmm1}.

Assume that  $(v,b,r,k,\lambda)=(q^{3}+1, q^{2}(q^{3}+1), q^{3}, q, q-1)$.
Then from Lemma \ref{AK2} we known that there are only two conjugacy classes of  subgroups  of  order $q(q-1)$ in $T$
and  $H_{1}=Q_{1}:K\leq T_\alpha$ and $H_{2}=Q_{2}:K\leq T_\alpha$ as representatives, respectively.

First, we consider the orbits of $H_{1}$. Let $\gamma\neq\alpha$ be the point fixed by $K$.
Since $K\leq H_{1}$, then $K_{\gamma}=K\leq (H_{1})_{\gamma}\leq T_{\alpha\gamma}=K,$
which implies ${(H_{1})}_{\gamma}=T_{\alpha\gamma}$ and so $|H_1:(H_{1})_{\gamma}|=|\gamma^{H_{1}}|=q$.
It is easy to see that  $|\delta^{H_{1}}|\neq q$ for any point $\delta\neq\alpha,\gamma$.
Therefore, $H_{1}$ has only one orbit of size $q$. Let $B_{1}=\gamma^{H_{1}}$.

Now we  show that $H_{1}= T_{B_{1}}$, which implies $|B_{1}^{T}|=b$.
Since $H_{1}\leq T_{B_{1}}$ and  $B_{1}=\gamma^{H_{1}}=\gamma^{T_{B_{1}}}$, then $|H_{1}:(H_{1})_{\gamma}|=|T_B:T_{\gamma B_{1}}|=q.$
If $K=(H_{1})_{\gamma}<T_{\gamma B_{1}}$,
then 3 divides $|T_{\gamma B_{1}}:T_{\delta\gamma B_{1}}|$ for any $\delta\in B_{{1}}\setminus \{\gamma\}$ by Lemma \ref{rpro}(3).
It follows that $3\mid (q-1)$,  a contradiction.
As a result, $K=(H_{1})_{\gamma}=T_{\gamma B_{1}}$ and so $H_{1}= T_{B_{1}}$. Let ${\mathcal{B}}_{1}:=B_{1}^{T}$.
Therefore $|{\mathcal{B}}_{1}|=|T:H_{1}|=b$. Let ${\mathcal{B}}_{1}$ be the set of blocks.

Finally, since $T$ is 2-transitive on  $\mathcal{P}$,   the number of blocks which incident with  two points is a constant.
Hence ${\mathcal{D}}_{1}=(\mathcal{P},{\mathcal{B}}_{1})$ is a  2-$(q^{3}+1,  q, q-1)$ design
admitting $T$ as a  flag transitive automorphism group by Lemma \ref{lemma2.1}.

In a  similar way,  we get the  design ${\mathcal{D}}_{2}$ satisfying all hypothesis when the subgroup is  $H_{2}=Q_{2}:K$.
Furthermore, since $H_{1}$  is not isomorphic to $H_{2}$, so $\mathcal{D}_{1}$ is not isomorphic to $\mathcal{D}_{2}$ by \cite[1.2.17]{Demb1968}.

Similarly , if   $(v,b,r,k,\lambda)=(q^{3}+1, q(q^{3}+1), q^{3}, q^{2}, q^2-1)$,
we can construct the  design with these parameters.
$\hfill\square$

\subsection{T is  the Suzuki group}

\begin{proposition} \label{pro2}
Suppose that $G$ and $\mathcal{D}$ satisfy hypothesis of Theorem \ref{theorem}.
If $T={}^{2}B_{2}(q)$ with $q=2^{2n+1}$,
then  $\cal D$ is a $2$-$(q^{2}+1, q, q-1)$ design with  $G_B=Z({Q}):K$ where $Q\in Syl_2(T)$ and $K=Z_{q-1}$.
\end{proposition}
{\bf Proof.}\, Suppose that $T={}^{2}B_{2}(q )$  with order $(q^{2}+1)q^{2}(q-1)$. Then $|G|=f(q^{2}+1)q^{2}(q-1)$ where $f$ divides $|Out(T)|$.
By \cite{XGFS} or \cite{Suzu}, the order of $G_\alpha$  is one of the following:
\begin{enumerate}
\item[\rm(1)] $fq^2(q-1)$;
\item[\rm(2)] $2f(q-1)$;
\item[\rm(3)] $4f(q\pm\sqrt{2q}+1)$;
\item[\rm(4)] $f(q_{0}^{2}+1)q_{0}^{2}(q_{0}-1)$ with $q_0^{\ell}=q$.
 \end{enumerate}

Since $|G|<|G_\alpha|^{3}$,  we first have that $|G_\alpha|\neq 2f(q-1)$.
If $|G_\alpha|=4f(q\pm\sqrt{2q}+1)$, from the inequality  $|G|<|G_\alpha|^{3}$, we get
$f(q^2+1)q^2(q-1)<(4f)^3(2q)^3$, and so $q^2+q+1\leq4^3f^2 2^3$. Since $f\leq |Out(T)|=e$ and $q=p^e$, hence $q+1<4^3 2^3$ and $q=2^7$, $2^5$ or $2^3$.
If $q=2^7$, then $ |G|=f2^{14}(2^{14}-1)(2^7-1)>f^34^3(2^7+2^4+1)^3=|G_\alpha|^3$  where $f=7$ or 1, a contradiction.
If $q=2^5$, then $v=198400$   or 325376 for $|G_\alpha|=4f(q+\sqrt{2q}+1)$ or  $4f(q-\sqrt{2q}+1)$ respectively.
By calculating $(|G_\alpha|, v-1)$, since $r$ divides $(|G_\alpha|, v-1)$, we know that $r$ is too small.
Similarly, we get  $q\neq 2^3$.

If $|G_\alpha|=f(q_{0}^{2}+1)q_{0}^{2}(q_{0}-1)$ with $q_0^{\ell}=q$, then the inequality $|G|<|G_\alpha||G_\alpha|_{p'}^{2}$ forces $m=3$.
So $v =(q_0^4-q_0^2+1)q_0^4(q_0^2+q_0 + 1)$.  Since $r$ divides $(|G_\alpha|_{p'}, v-1)$, then $r \leq  |G_\alpha|_{p'} \leq  fq_0^3 <q_0^{9/2}$.
From $v<r^2$,  we get $(q_0^4-q_0^2+1)q_0^4(q_0^2+q_0 + 1)<r^2<q_0^{9}$, which is impossible.

Now assume that $|G_\alpha|=fq^2(q-1)$.  Then $v=q^{2}+1$ and $T$ is 2-transitive by \cite[p.250]{pg}.
Hence, $T$ is flag-transitive   by Lemma \ref{do}.
Similarly, we have   $|T_{B}|=\frac{|T|}{b}=\frac{k(k-1)(q-1)}{\lambda}.$
Let $M$ be the maximal subgroup of $T$ such that $T_B\leq M$  as in Lemma \ref{lmmm1}. The fact that  $|T_{B}|$ divides $ |M|$ implies that $|M|=q^{2}(q-1)$
and $k(k-1)$ divides $\lambda q^{2}$.
 Similar to the proof of  Lemma \ref{lmmm1}, we have $T_{\gamma B}=T_{\alpha\gamma}$  with the order $q-1$.
Furthermore, we get $$(v,b,r,k,\lambda)=(q^{2}+1, q(q^{2}+1), q^{2}, q, q-1).$$

Next we prove  the existence of the design  with above parameters by Lemma \ref{lemma2.1}.
Firstly, from Lemma \ref{AK4} we know  that  the Suzuki group has a unique conjugacy class of subgroups of order $q(q-1)$, let $H:=Z(Q):K\leq T_\alpha$ as the representative.

Note that  $K$ is the stabilizers of two points in $^{2}B_{2}(q)$ by \cite[p.187]{BH1982}. Let $\gamma\neq\alpha$ be the point fixed by $K$ and $B=\gamma^{H}$.
Then similar as the proof of Proposition \ref{pro1} we get  that $B$ is the only  $H$-orbit of length $q$ and $H=T_B$.
Let $\mathcal{B}=B^{T}$ be the set of blocks.
Finally, since $T$ is 2-transitive on  $\mathcal{P}$,   the number of blocks which incident with  two points is a constant.
Hence $\mathcal{D}=(\mathcal{P}, \mathcal{B})$ is a  2-$(q^{2}+1,  q, q-1)$ design
admitting $T$ be a  flag transitive automorphism group by Lemma \ref{lemma2.1}.
$\hfill\square$

\subsection{T is one of   the remaining   families }
In this subsection, let  $$ \mathcal{T}=\{{^{2}F_{4}(q)},  {^{3}D_{4}(q)}, G_{2}(q), F_4(q), E_6^\epsilon(q), E_7(q), E_8(q)\},$$
 we will prove that there are no new design arise when $T\in\mathcal{T}$.

First, we show that $G_{\alpha}$ cannot be   a parabolic subgroup of $G$  for any $T\in\mathcal{T}$.
\begin{lemma} \label{pro1.1}
Suppose that $G$ and $\mathcal{D}$ satisfy hypothesis of Theorem \ref{theorem}. If  $T\in  \mathcal{T}$, then $G_{\alpha}$ cannot be   a parabolic subgroup of $G$.
\end{lemma}
{\bf Proof.}\, By  Lemma \ref{lem:power}, for all cases that $T\in  \mathcal{T}\setminus E_{6}(q)$, there is a unique subdegree which is a power of $p$,  so  $r$ is a  power of $p$ by Lemma \ref{lemma1.2}(4). We can easily check that  $r$ is too small and   the condition $r^{2}>v$ cannot be satisfied.
Now, assume that  $T = E_{6}(q)$.
If $G$ contains a graph automorphism or $G_{\alpha}\cap T$ is $P_{2}$ or $P_{4}$,
then there is also a unique subdegree which  is a power of $p$ and so $r$ is too small again.
If $G_{\alpha}\cap T$ is $P_{3}$ with type $A_{1}A_{4}$,
then $$v=\frac{(q^{3}+1)(q^{4}+1)(q^{9}-1)(q^{6}+1)(q^{4}+q^{2}+1)}{(q-1)}.$$
Since $r$ divides $(|G_{\alpha}|,v-1)$, we have
$r\mid eq(q-1)^{5}(q^{5}-1)$ and  so $r$ is too small to satisfy $r^{2}>v$.
If $G_{\alpha}\cap T$ is $P_{1}$ with type $D_{5}$,
then $$v=\frac{(q^{8}+q^{4}+1)(q^{9}-1)}{q-1}.$$
From $\cite{LieM1986}$, we know that there exists two non-trivial subdegrees:
$$d=\frac{q(q^{3}+1)(q^{8}-1)}{(q-1)}~~~\mathrm{and}~~d'=\frac{q^{8}(q^{4}+1)(q^{5}-1)}{(q-1)}.$$
Since  $(d,d')=q(q^{4}+1)$,  we have $r\mid q(q^{4}+1)$ by Lemma \ref{lemma1.2}(4), which contradicts  with $r^{2}>v$.$\hfill\square$

Let $\mathcal{T}_{1}=\{F_4(q), E_6^\epsilon(q), E_7(q), E_8(q)\}.$

\begin{lemma} \label{maxrk}
Suppose that  $G$ and $\mathcal{D}$ satisfy the hypothesis of Theorem \ref{theorem}. If   $T\in \mathcal{T}_{1}$ and $G_{\alpha}$ is non-parabolic, then    $G_\alpha$ cannot be a maximal subgroup of maximal rank.
\end{lemma}
{\bf Proof.}\, If $G_\alpha$ is non-parabolic and of maximal rank,  then for any  $T\in \mathcal{T}_{1}$,
we have a complete list of $T_\alpha:=G_\alpha \cap T$  in \cite[Tables 5.1-5.2]{LSS92}.
All subgroups   in \cite[Table 5.2]{LSS92} and some cases in  \cite[Table 5.1]{LSS92} can be ruled out by the inequality $|T|<|Out(T)|^2|T_\alpha|{|T_\alpha|_{p'}^{2}}$.
Since $r$ divides $(|G_\alpha|, v-1)$, for  the remaining cases we have that $r^2<v$, a contradiction.

For example, if $T=F_4(q)$ with order $q^{24}(q^{2}-1)(q^{6}-1)(q^{8}-1)(q^{12}-1)$.  Then $T_\alpha$ is one of the following:
(1) $2.(L_2(q)\times PSp_6(q)).2$  ($q$ odd);
(2) $d.\Omega_9(q)$; (3) $d^{2}.P\Omega^{+}_{8}(q).S_{3}$; (4) $^{3}D_{4}(q).3$; (5) $Sp_{4}(q^{2}).2$ ($q$ even);
(6) $(Sp_{4}(q)\times Sp_{4}(q)).2$($q$ even); (7) $h.(L^{\epsilon}_{3}(q)\times L^{\epsilon}_{3}(q)).h.2$,
with $d=(2,q-1)$ and $h=(3,q-\epsilon)$.

If $T_{\alpha}=2.(L_2(q)\times PSp_6(q)).2$ with $q$ odd,  then$$|T_{\alpha}|=q^{10}(q^{2}-1)^{2}(q^{4}-1)(q^{6}-1)~~\mathrm{and}~~
v=q^{14}(q^{4}+1)(q^{4}+q^{2}+1)(q^{6}+1).$$
Since $(q^{2}+1)\mid v$ and $(q^{4}+q^{2}+1)\mid v$, then $(|G_\alpha|, v-1)\mid |Out(T)|(q^{2}-1)^{4}$ and so $r^{2}<q^9<v$,  a contradiction.

If $T_{\alpha}=2.P\Omega_{9}(q)$ with $q$ odd, then
$$|T_{\alpha}|=q^{16}(q^{2}-1)(q^{4}-1)(q^{6}-1)(q^{8}-1)~~\mathrm{and}~~v=q^{8}(q^{8}+q^{4}+1).$$
Since  $q\mid v$, $(q^{4}+q^{2}+1)\mid v$, $v-1\equiv 2 {\pmod {q^{4}-1}}$,  we get $r$ divides $2^{4}|Out(T)|(q^{4}+1)$ and so $r^{2}<v$,  a contradiction.

 Cases (3)-(6) can be  ruled out similarly, and Case (7) cannot occur because of $|T|<|Out(T)|^2|T_\alpha|{|T_\alpha|_{p'}^{2}}$.       $\hfill\square$

\begin{lemma} \label{notsimple}
Suppose that  $G$ and $\mathcal{D}$ satisfy the hypothesis of Theorem \ref{theorem}. If   $T\in \mathcal{T}_{1}$ and $G_{\alpha}$ is non-parabolic,
then  $T_0=Soc(G_\alpha\cap T)$ is  simple and  $T_0=T_0(q_0)\in Lie(p)$.
\end{lemma}
{\bf Proof.}\, Assume that  $T_0=Soc(G_\alpha\cap T)$ is not  simple. Then by Lemma \ref{theorem 10} and Lemma \ref{maxrk}, one of the following holds:
\begin{enumerate}
\item[\rm(1)] \, $G_\alpha= N_G(E)$, where $E$ is an elementary abelian group given in \cite[Theorem 1(II)]{CLSS92};
\item[\rm(2)] \, $T=E_8(q)$ with $p>5$, and $T_0$ is either $A_5\times A_6$ or $A_5\times L_2(q)$;
\item[\rm(3)] \, $T_0$ is as in Table \ref{Table 1}.
\end{enumerate}

From \cite[Theorem 1(II)]{CLSS92}, we check  that  all subgroups in Case (1) are local and  too small to satisfy $|T|<|Out(T)|^2|T_\alpha|{|T_\alpha|_{p'}^{2}}$.


The order of subgroup in Case (2) is too small.

For Case (3),  since $G_\alpha$ is not simple and  not local by \cite[Theorem 1]{CLSS92},  $G_\alpha$ is of maximal rank by \cite[p.346]{JS2001}, which has already been ruled out in Case (1).
Therefore,  $T_0$ is simple.

Now assume that  $T_0=T_0(q_0)\not\in Lie(p)$. Then for all $T$, we find the
possibilities of $T_0$ in   \cite[Table 1]{LS99}.  Some cases can be ruled out by the inequality $|T|<|Out(T)|^2|T_\alpha|{|T_\alpha|_{p'}^{2}}$.
In each of the remaining cases,
since $r$ must divides $(|G_\alpha|, v-1)$, $r$ is too small to satisfy $v<r^2$.
For example, assume that  $T= F_4(q)$. If  $T_0\not\in Lie(p)$, then  according to   \cite[Table 1]{LS99}, it is one of the following:
 $A_{5-10}$, $L_2(7)$, $L_2(8)$, $L_2(13)$, $L_2(17)$, $L_2(25)$, $L_2(27)$, $L_3(3)$, $U_3(3)$, $U_4(2)$,
 $Sp_6(2)$, $\Omega_8^+(2)$, $ {^3D_4(2)}$, $J_2$, $J_2$, $A_{11}(p=11)$, $L_3(4)(p=3)$,
 $L_4(3)(p=2)$, ${^2B_2(8)}(p=5)$, $M_{11}(p=11)$.
The  possibilities of $T_0$ such that $|G|<|G_\alpha|^3$ are $A_9(q=2)$, $A_{10}(q=2)$, $Sp_6(2)(q=2)$, $\Omega_8^+(2)(q=2,3)$, ${^3D_4(2)}(q=2,3)$, $J_2(q=2)$, $L_{4}(3)(q=2)$.
However, since $r\mid(|G_\alpha|, v-1)$, we have $r^2<v$ for all these cases, which is a contradiction.$\hfill\square$

\begin{lemma} \label{simple}
Suppose that  $G$ and $\mathcal{D}$ satisfy the hypothesis of Theorem \ref{theorem}. If $T_0=T_0(q_0)$ is a simple group of Lie type  and  $G_{\alpha}$ is non-parabolic,  then  $T\not\in \mathcal{T}_{1}.$
\end{lemma}
{\bf Proof.}\, First assume that $T= F_4(q)$.   If ${\rm rank}(T_0)> \frac{1}{2}{\rm rank}(T)$,
then by Lemma \ref{theorem 11} and Lemma \ref{maxrk}, the only possible cases of $G_\alpha\cap T$ satisfying $|G|<|G_\alpha|^3$
are $F_4(q^{\frac{1}{2}})$ and $F_4(q^{\frac{1}{3}})$ when $q_{0}>2$.
If $G_\alpha\cap T=F_4(q^{\frac{1}{2}})$, then $v=q^{12}(q^6+1)(q^4+1)(q^3+1)(q+1)>q^{26}.$
Since $q$, $q+1$, $q^{2}+1$ and $q^{3}+1$ are factors of $v$,
then $r\mid 2e(q-1)^{2}(q^{3}-1)^{2}$ by $r\mid (|G_\alpha|, v-1)$, which implies that $r^{2}<v$, a contradiction.
If $G_\alpha\cap T=F_4(q^{\frac{1}{3}})$, since $p\mid v$, then $r$ divides $|G_\alpha|_{p^{'}}$, which also implies $r^2<v$.
When  $q_{0}=2$,  the  subgroups  $T_0(2)$ with ${\rm rank}(T_0)> \frac{1}{2}{\rm rank}(T)$ that satisfy   $|G|<|G_\alpha|^3$ are $A_4^\epsilon(2)$,  $B_3(2)$, $B_4(2)$, $C_3(2)$, $C_4(2)$ or $D_4^\epsilon(2)$.
But in each case, $r\mid  (|G_\alpha|, v-1)$ forces $r^2<v$, a contradiction.
If ${\rm rank}(T_0)\leq \frac{1}{2}{\rm rank}(T)$, then  from  Lemma \ref{theorem 12}, we have $|G_\alpha|<4 q^{20}\log_pq$.
Looking at the orders  of   groups of Lie type, we see that if $|G_\alpha|<4 q^{20}\log_pq$, then  $|G_\alpha|_{p'}<q^{12}$, and so $|G_\alpha||G_\alpha|_{p'}^2<|G|$, contrary to  Lemma \ref{lifanggen}.

For $T= E_6^\epsilon(q)$, if ${\rm rank}(T_0)>\frac{1}{2} {\rm rank}(T)$, then when  $q_0>2$, by Lemma \ref{theorem 11} the only possibilities are $E^{\epsilon}_6(q^{\frac{1}{2}})$, $E^{\epsilon}_6(q^{\frac{1}{3}})$, $C_4(q)$ and $F_4(q)$.
In all these cases $r$ are too small.
When $q_0=2$, the possibilities $T_0(2)$ satisfying $|G|<|G_\alpha|^3$ with order dividing
$|E_6^\epsilon(2)|$ are  $A^\epsilon_5(2)$, $B_4(2)$, $C_4(2)$, $D_4^\epsilon(2)$  and $D_5^\epsilon(2)$. However, since $r\mid  (|G_\alpha|, v-1)$, for all these cases we obtain $r^2<v$, a contradiction.
If ${\rm rank}(T_0)\leq \frac{1}{2}{\rm rank}(T)$, then  from  Lemma \ref{theorem 12}, we have $|G_\alpha|<4 q^{28} \log_pq$.
By further check the orders of  groups of Lie type, we see that   $|G_\alpha|_{p'}<q^{17}$, and so $|G_\alpha||G_\alpha|_{p'}^2<|G|$, a contradiction.

Assume that $T= E_7(q)$. If $rank(T_0)\leq\frac{1}{2}rank(T)$, then by Lemma \ref{theorem 12}
$|G_\alpha|^3\leq|G|$, a contradiction.
If ${\rm rank}(T_0)> \frac{1}{2}{\rm rank}(T)$,
then when $q_{0}>2$, B by Lemma \ref{theorem 11}, the only cases $T\cap G_\alpha$ satisfying $|G|<|G_\alpha|^3$ are $G_\alpha\cap
T=E_7(q^{\frac{1}{s}})$, where $s=2$ or $3$.  But in all cases we have $r^2<v$.
If $q_0=2$, then the possible subgroups such that $|G|<|G_\alpha |^3$ with order dividing $|E_7(2)|$
are $A_6^\epsilon(2)$, $A_7^\epsilon(2)$, $B_5(2)$, $C_5(2)$, $D_5^\epsilon(2)$ and $D_6^\epsilon(2)$.
However in all of these cases, since $r\mid (|G_\alpha|, v-1)$ we have  $r^2<v$, a contradiction.

Assume that $T= E_8(q)$. If ${\rm rank}(T_0)\leq \frac{1}{2}{\rm rank}(T)$, then by Lemma \ref{theorem 12}
we get $|G_\alpha|^3<|G|$, a contradiction. Therefore, ${\rm rank}(T_0)>\frac{1}{2} {\rm rank}(T)$.
If $q_0>2$, then Lemma \ref{theorem 11} implies  $G_\alpha\cap T=E_8(q^{\frac{1}{s}})$, with $s=2$ or $3$.
However in both cases we get a small $r$ with $r^2<v$, a contradiction.
If $q_0=2$, then ${\rm rank}(X_0)\geq 5$. All subgroups satisfying $|G_\alpha|^3>|G|$ are $A_8^\epsilon(2)$, $B_7(2)$, $B_8(2)$, $C_7(2)$, $C_8(2)$, $D_8^\epsilon(2)$ and $D_7^\epsilon(2)$. But for all these cases  we have  $r^2<v$.$\hfill\square$

\begin{lemma}\label{g2}
If $T=G_{2}(q)$ with $q=p^e>2$, then $G_{\alpha}$ cannot be  a non-parabolic maximal subgroup of $G$.
\end{lemma}
{\bf Proof.}\, Suppose that $T=G_2(q)$ with $q>2$ since $G_{2}(2)'=PSU_{3}(3)$.
All maximal subgroups of $G$ can be found in \cite{PBK1} for odd $q$ and in \cite{BNC} for even $q$.

Assume that $G_{\alpha}$  be  a non-parabolic maximal subgroup of $G$. First we deal with the case where   $G_\alpha\cap T=SL_3^\epsilon(q).2$ with $\epsilon=\pm$.
Then we have  $v=\frac{1}{2}q^{3}(q^{3}+\epsilon1).$
By Lemma \ref{lemma1.2} and \cite[Section 8]{JS2001} we conclude  that
$r$ divides $\frac{(q^{3}-\epsilon1)}{2}$ for odd $q$ (cf. \cite[Section 4, Case 1, $i=1$]{JS2001})
and $r$ divides $(q^{3}-\epsilon1)$ for even $q$ (cf. \cite[Section 3, Case 8]{JS2001}).
The case that $q$ odd is ruled out by $v<r^{2}$.  If $q$ is even, then $r=q^{3}-\epsilon1$.
This, together with $k<r$,  implies   $k-1=\lambda\frac{q^{3}+\epsilon2}{2}$, and so $\lambda=1$ or $\lambda=2$.
From the result of \cite{JS2001} we known that $\lambda\neq1$. If $\lambda=2$, then since $k<r$, we have $\epsilon=-$.
It follows that $k=q^{3}-1$ and $r=q^{3}+1$. This is impossible by  Lemma \ref{lemma115} and \cite[Theorem 1]{EOL}.

Now, if $G_\alpha\cap T={^2G_2(q)}$ with $q=3^{2n+1}\geq27$, then
$v=q^{3}(q+1)(q^{3}-1).$
Note that $q\mid v$ and $(q^2-1,v-1)=1$, we have  $(|G_{\alpha}|, v-1)\mid e(q^{2}-q +1)$,
and it follows that $r^{2}<v$,  a contradiction.

The cases  that $G_\alpha\cap T$ is $G_2(q_0)$ or $(SL_2(q)\circ SL_2(q))\cdot 2$ can be ruled out similarly.

Using  the inequality $|G|<|G_\alpha|^3$ and the fact that $r$ divides $(|G_{\alpha}|,v-1)$,
we find  $r$ too small to satisfy $ r^2>v$ for every other    maximal subgroup.
$\hfill\square$

\begin{lemma}\label{f24}
If $T={^{2}F_{4}(q)}$, then $G_{\alpha}$ cannot be  a non-parabolic maximal subgroup.
\end{lemma}
{\bf Proof.}\, Let  $T={^{2}F_{4}(q)}$ and $G_{\alpha}$   be  a non-parabolic maximal subgroup of $G$.
Then from the list of the maximal subgroups of $G$ in \cite{MG},
there are no   subgroups $G_\alpha $ satisfying   $|G|<|G_\alpha|{|G_\alpha|_{p'}^{2}}$,
except for the case $q=2$. For the case $q=2$, $G_\alpha\cap T$ is $ L_3(3).2$ or $ L_2(25)$.
However in each case, since $r$ divides $(|G_\alpha|, v-1)$,  and so $r$ is too small.$\hfill\square$

\begin{lemma}\label{d4}
If $T={^{3}D_{4}(q)}$, then $G_{\alpha}$ cannot be  a non-parabolic maximal subgroup.
\end{lemma}
{\bf Proof.}\,
If  $T={^{3}D_{4}(q)}$ and $G_{\alpha}$  is  a non-parabolic maximal subgroup of $G$, then all possibilities of $G_\alpha\cap T$ are listed in \cite{PBK2}.
However, for all cases, the fact that $r$ divide $(|G_{\alpha}|, v-1)$ give a small $r$ which cannot satisfy the condition  $v<r^{2}$.
For example, if  $G_\alpha\cap T$ is $G_{2}(q)$ of order $q^{6}({q}^2-1)({q}^6-1)$, then $v={q}^{6}({q}^8+q^4+1)$.
Since $q\mid v $ and $(q^{4}+q^{2}+1)\mid v$, then $r\mid3e(q^{2}-1)^{2}$, which contradicts with $v<r^{2}$. $\hfill\square$
\begin{lemma} \label{pro2.1}
Suppose that  $G$ and $\mathcal{D}$ satisfy the hypothesis of Theorem \ref{theorem}. If the socle $T\in \mathcal{T}$,
then $G_{\alpha}$ cannot be a non-parabolic maximal subgroup.
\end{lemma}

{\bf Proof.}\, It is follows from Lemmas \ref{maxrk}$-$\ref{d4}. $\hfill\square$

{\bf Proof of Theorem  \ref{theorem}.}\,  Now Theorem \ref{theorem} is an immediate consequence of Propositions \ref{pro1}-\ref{pro2} and of Lemmas \ref{pro1.1} and \ref{pro2.1}.

\end{document}